\documentclass[12pt,titlepage]{article}

\usepackage{amssymb}
\usepackage{exscale}
\usepackage{url}
\usepackage[amsthm]{ntheorem}

\input epsf
\newcommand{\gproof}{{\noindent \bf Proof. }}
\theoremstyle{plain}
\newtheorem{thm}{Theorem}
\newtheorem{defi}{Definition}

\newtheorem{cor}[thm]{Corollary}
\newtheorem{prop}[thm]{Proposition}

\newtheorem*{tquestion}{Question}
\makeatletter
\newtheoremstyle{named}
	{\item[\hskip\labelsep\theorem@headerfont##1]}%
	{\item[\hskip\labelsep\theorem@headerfont##1\ ##3]}
\newtheoremstyle{nonumbernamed}
	{\item[\hskip\labelsep\theorem@headerfont##1]}%
	{\item[\hskip\labelsep\theorem@headerfont##1\ ##3]}
\newtheoremstyle{nnamed}
	{\item[\hskip\labelsep\theorem@headerfont##1]}%
	{\item[\hskip\labelsep\theorem@headerfont##3]}
\newtheoremstyle{nonumbernnamed}
	{\item[\hskip\labelsep\theorem@headerfont##1]}%
	{\item[\hskip\labelsep\theorem@headerfont##3]}
\theorembodyfont{\itshape}
\theoremseparator{.}
\makeatother
\theoremstyle{named}
\newtheorem*{namedthm}{Theorem}
\theoremstyle{nnamed}
\newtheorem*{nnamedthm}{Theorem}
\theoremstyle{remark}
\newtheorem{rem}{Remark}
\theoremstyle{plain}

\def\2{\mathbb Z_2}
\def\S{\mathbb S}
\def\ind{{\rm ind}}
\def\coind{{\rm coind}}

\def\KG{{\rm KG}}
\def\SG{{\rm SG}}

\title{On topological relaxations of chromatic conjectures}
\author{
      {\bf G\'abor Simonyi}\thanks{Research partially supported by the
Hungarian Foundation for Scientific Research Grant (OTKA) Nos.\ 
K76088 and NK78439}\\
 Alfr\'ed R\'enyi Institute of Mathematics\\
 Hungarian Academy of Sciences\\
 {\tt simonyi@renyi.hu}
\and
   {\bf Ambrus Zsb\'an}\thanks{Research partially supported by the Hungarian
     Foundation for Scientific Research Grant and by the National
Office for Research and Technology (Grant number OTKA 67651).}\\ 
 Department of Computer Science and Information Theory\\
 Budapest University of Technology and Economics \\
 {\tt ambrus@math.bme.hu}}
\date{}

\begin{document}
\maketitle

\begin{abstract}
There are several famous unsolved conjectures about the chromatic number that
were relaxed and already proven to hold for the fractional chromatic number. We
discuss similar relaxations for the topological lower bound(s) of the chromatic
number. In particular, we prove that such a relaxed version is true for the
Behzad-Vizing conjecture and also discuss the conjectures of Hedetniemi and of
Hadwiger from this point of view. For the latter, a similar statement was
already proven in \cite{ST}, 
our main concern here is that the so-called odd Hadwiger conjecture looks much
more difficult in this respect. We prove that the statement of the
odd Hadwiger conjecture holds for large enough Kneser graphs and Schrijver
graphs of any fixed chromatic number.  
\end{abstract}

\section{Introduction}

There are several hard conjectures about the chromatic number that
are still open, while their fractional relaxation is solved, i.e., a 
similar, but weaker statement is proven for the fractional chromatic
number in place of the chromatic number. (For the definition and basic facts
about the fractional chromatic number, cf.\ \cite{SchU}.) Examples include the
Behzad-Vizing conjecture 
\cite{KR}, the Erd\H{o}s-Faber-Lov\'asz conjecture \cite{KahnSeym},
Hedetniemi's conjecture \cite{THed2}, a relaxed version of Hadwiger's
conjecture \cite{RS}, as well as a similarly relaxed version of the so-called
odd Hadwiger conjecture \cite{KKR}. (In some of these cases the proven
fractional version has an approximative form, nevertheless, it is a statement
not known to hold for the chromatic number.)  

There are not very many examples of graphs with a large gap between their
chromatic and fractional chromatic numbers. To determine the chromatic number
of such a graph is usually difficult because no lower bound that also bounds
the fractional chromatic number from below can give a tight result. The
primary example for such a graph family is that of Kneser graphs. The value of
their chromatic number was conjectured by Kneser \cite{Kne} in 1955 and proved
by Lov\'asz \cite{LLKn} in 1978 thereby developing the topological method for
estimating the chromatic number. This method was later successfully applied to
other graphs, e.g., generalized Mycielski graphs, cf.\ \cite{Stieb, GyJS}, see 
\cite{Mat} for a thorough survey on the later developments. 

The above suggests that one could gain further supporting evidence for the
above conjectures if one could prove that the topological lower bound for the
chromatic number, considered as a graph parameter for its own sake, also
satisfies the above statements if it is put in place of the chromatic number. 
(In fact, ``the topological lower bound'' is not a well-defined term, as there
are more than one such bounds, for further details see the next section.) 

Such a result already appears in the last section of \cite{ST} concerning
a relaxation of Hadwiger's conjecture. Further impetus for such studies was
given to us by a conversation of the first author with Claude Tardif and
G\'abor Tardos at the first Canadam Conference in summer 2007 where the idea
of considering the 
topological lower bound(s) as a graph parameter was made more explicit. In
particular, Tardif asked, whether a result of the above type would be possible
concerning the Erd\H{o}s-Faber-Lov\'asz conjecture. Though we were not able to
make progress in this particular question, we will prove in this paper a
similar result about the Behzad-Vizing conjecture and elaborate about some of
the others. 

The paper is organized as follows. Section~\ref{top} contains some basic facts
about 
the topological method. In Section~\ref{secBH} we prove our result concerning
the 
Behzad-Vizing conjecture. In Section~\ref{secHedet} we give a simple
topological analog of Hedetniemi's conjecture. In Sections~\ref{secHad} and 
\ref{secoddHad} we discuss Hadwiger's conjecture and the odd Hadwiger
conjecture. In the latter section we prove that the odd Hadwiger conjecture
holds for some of the graphs for which the topological method gives a tight
bound on the chromatic number, in particular, large enough Kneser graphs,
Schrijver graphs, and generalized
Mycielski graphs.

\section{About the topological bound(s) on the chromatic number} \label{top}

There are several formally different topological lower bounds on the chromatic
number that are all closely related to each other. As we will only use
combinatorial consequences of the situation when these parameters achieve
certain values rather than using them directly, we will not give full
definitions of these bounds. (Most importantly, we are not defining the
topological notions used. They can be looked up in several of the references
and though they give important background, familiarity with these notions, or
even knowing them, is not essential for understanding this paper.) 
Instead we only hint the definitions and give references for detailed
treatments, while list those statements that are to be used in this paper. 

The idea behind all versions of the topological lower bound of the chromatic
number is to associate a topological space to the graph and use its topological
invariants for bounding the chromatic number. Originally Lov\'asz \cite{LLKn} 
used the connectivity of the associated topological space defined via a
simplicial complex, called the neighborhood complex, and showed that this
parameter is less than the chromatic number of the graph by at least $3$.  

Other variants of the same idea appeared over the years that use $\2$-spaces
defined by certain box complexes. (In fact, $\2$-maps and the Borsuk-Ulam
theorem are also key in Lov\'asz' original proof, the difference is only the
more direct use of $\2$-spaces in these later variants.) 
For a variety of box complexes, see
\cite{MZ}. 
One of the most basic box complexes, $B(G)$, associated to graph
$G$, is a simplicial complex that has $V(G)\times\{1\}\cup V(G)\times\{2\}$ as
its vertex set, and a subset of vertices forms a simplex in it iff it has the
form $A\uplus B:=A\times\{1\}\cup B\times\{2\}$, the induced subgraph of $G$
on $A\cup B\subseteq V(G)$ contains a complete bipartite graph with color classes
$A$ and $B$, and in case $A$ (or $B$) is empty, we have that the vertices in
$B$ (resp. $A$) have at least one common neighbor. (In other words, simplices
of the form $A\uplus\emptyset$ and $\emptyset\uplus B$ are contained only when
they should be by the hereditary nature of simplicial complexes.) 
The $\2$-space evolving from
this simplicial complex is the topological space given by its geometric
realization equipped with the $\2$-map generated by the simplicial map $\nu:
A\uplus B\mapsto B\uplus A$. The most important property of this construction
is that whenever $G$ and $H$ are two graphs such that there exists a
homomorphism, i.e., an edge preserving map of the vertices, from $G$ to $H$,
then there is also a simplicial $\2$-map from $B(G)$ to $B(H)$. It is not
hard to show that $B(K_n)$ is homotopy equivalent to the sphere
$\S^{n-2}$. (As 
a $\2$-space the sphere $\S^h$ is considered to be equipped with the antipodal
map as the $\2$-map.) One can define for any $\2$-space $T=(T,\nu)$ its
$\2$-index $\ind(T)$ as the smallest dimension $h$ for which a $\2$-map
(that is, one respecting the involution $\nu$) exists from $T$ to $\S^h$. The
celebrated 
Borsuk-Ulam theorem (cf.\ e.g.\ \cite{Mat}) states in one of its standard forms,
that no antipodal map exists from $\S^h$ to $\S^{h'}$ if $h'<h$. Putting all
this together, and using the fact that a proper coloring with $m$ colors is
nothing but a homomorphism to $K_m$, one obtains that $\chi(G)\ge
\ind(B(G))+2$ should always hold.  
 
A somewhat different box complex $B_0(G)$ can be
defined by simply dropping the extra condition about common neighbors for
the containment of simplices having the form $A\uplus\emptyset$ or
$\emptyset\uplus B$. (Thus $V(G)\uplus\emptyset$, $\emptyset\uplus V(G)$ and
all their subsets are simplices in $B_0(G)$.)  
Csorba \cite{CsCCA} proved that  
$B_0(G)$ is $\2$-homotopy equivalent to the suspension of
$B(G)$, 
cf.\ also \cite{MZ} for this and other relations between various box
complexes.     
The latter fact implies $\ind(B_0(G))\le\ind(B(G))+1$ thus we have by
the foregoing that the inequality $\chi(G)\ge \ind(B_0(G))+1$ holds, too.

By the above mentioned form of the Borsuk-Ulam theorem, 
the $\2$-index of a $\2$-space $T$ is bounded from below by the $\2$-coindex,
$\coind(T)$ which is defined as the largest dimension $h$ for which a $\2$-map
exists {\em from} the sphere $\S^h$ to $T$. By the suspension relationship we
have $\coind(B_0(G))\ge \coind(B(G))+1$. Thus the $\2$-index and $\2$-coindex
of the two box complexes we discussed give the following chain of lower bounds
on the chromatic number: 
$$\chi(G)\ge
\ind(B(G))+2\ge\ind(B_0(G))+1\ge\coind(B_0(G))+1\ge\coind(B(G))+2.$$ 
For a more thorough introduction to these notions we refer to \cite{Mat} or
\cite{ST}. 

Seeing the four lower bounds on $\chi(G)$ in the above chain of inequalities
one may ask why we do not keep only the strongest one and drop the rest. The
reason is that if a weaker lower bound of the above gives the same value as
one of the stronger ones, that may have stronger graph theoretic consequences
compared to the situation when there is a gap between the two bounds. An
example of this phenomenon is demonstrated in \cite{STV}. 

We will use the following results of earlier papers that give graph theoretic
consequences of the property that one of the above lower bounds attain a
certain value. 

The first such theorem we need involves the strongest of the above bounds. It
is proven by Csorba, Lange, Schurr, and Wa{\ss}mer in \cite{CsLSW} where it is
called the $K_{\ell,m}$-theorem. 

\begin{nnamedthm}[{\mathversion{normal}$K_{\ell,m}$-theorem}] \textup{(\cite{CsLSW})}
{\it If $G$ is a graph satisfying $\ind(B(G))+2\ge t$, then for every possible
$\ell,m\in\mathbb N$ with $\ell+m=t$, the complete bipartite graph
$K_{\ell,m}$ appears as a subgraph of $G$.} 
\end{nnamedthm}

The following result, that was named Zig-zag Theorem in \cite{ST}, involves
the third of the above bounds. 

\begin{nnamedthm}[Zig-zag Theorem] \textup{(\cite{ST}, cf.\ also \cite{kyfan2})}
{\it If $G$ is a graph satisfying $\coind(B_0(G))+2\ge t$, then the following
  holds for every proper coloring $c: V(G)\to\mathbb N$. $G$ contains a
  $K_{\lceil t/2\rceil,\lfloor{t/2}\rfloor}$ subgraph all $t$ vertices of
  which receive a different color by $c$. Furthermore, these $t$ colors, if
  considered in their natural order as natural numbers, appear alternately on
  the two sides of the given $K_{\lceil t/2\rceil,\lfloor{t/2}\rfloor}$
  subgraph.}
\end{nnamedthm}

Note that the number of colors used for the coloring in the Zig-zag Theorem
may be much more than 
$\chi(G)$. In case $\chi(G)=t=\coind(B_0(G))+1$ a colorful version of the
$K_{\ell,m}$-theorem is proven in \cite{STcol}.  
  
\medskip

We quote another result from \cite{ST} that gives a characterization of those
graphs for which the fourth of the above lower bounds is above a certain
value. This characterization needs the notion of Borsuk graphs defined by
Erd\H{o}s and Hajnal \cite{EH}. 

\begin{defi} \label{defi:Bogr} {\rm (\cite{EH})}
The Borsuk graph $B(n,\alpha)$ of parameters $n$ and $0<\alpha<2$ is an 
infinite graph whose vertices are the points of the unit sphere $\S^{n-1}$ in
${\mathbb R}^n$ and whose edges connect its pairs of vertices
with distance at least $\alpha$. 
\end{defi}

\begin{namedthm}[B]
\textup{(Lemma 4.4 in \cite{ST})}
{\it A finite graph $G$ satisfies $\coind(B(G))\ge n-1$ if and only if
there is a graph homomorphism from $B(n,\alpha)$ to $G$ for some $\alpha<2$.}
\end{namedthm}

Note that $\coind(B(G))\ge t-2$ implies $\chi(G)\ge t$ 
and $\chi(B(t-1,\alpha))=t$ for large enough
$\alpha<2$ is equivalent to the Borsuk-Ulam theorem, cf. \cite{LLgomb}.   
We remark that graphs satisfying
$\coind(B(G))\ge t-2$ are called strongly topologically $t$-chromatic in
\cite{ST} as opposed to topologically $t$-chromatic graphs defined by
satisfying $\coind(B_0(G))\ge t-1$.

\section{On the Behzad-Vizing conjecture} \label{secBH}

The Behzad-Vizing conjecture states that one can always color the vertices
and the edges of a simple graph $G$ with at most $\Delta(G)+2$ colors in such
a way that neither adjacent vertices nor edges with a common endvertex get the
same color, furthermore, no edge is colored the same as one of its
endpoints. Here 
$\Delta(G)$ denotes the maximum degree of $G$. The minimal number of colors
needed for such a coloring is called the {\em total chromatic number} and
is often denoted by $\chi''(G)$. It is simply the chromatic number of $T(G)$,
the {\em total graph} of $G$ defined by $$V(T(G))=V(G)\cup E(G)$$ and
$$E(T(G))=\{\{a,b\}: a,b\in V(G), \{a,b\}\in E(G)\ {\rm or}$$ $$a\in
V(G), b\in E(G), a\in b\ {\rm or}\ a,b\in E(G), a\cap b\neq\emptyset\}.$$ 

This problem is open for more than forty years.  Its original appearance seems
to be independently \cite{Beh, Beh2} and \cite{Viz}, see also 
\cite{OPG, IK}. It was solved for $\Delta(G)=3$ by 
Rosenfeld \cite{Ros} and Vijayaditya \cite{Vij} (it is trivial for
$\Delta(G)\le 2$)  
and for $\Delta(G)=4$ and $5$ by 
Kostochka \cite{Kost, Kost2, Kost3}. 
The fractional chromatic number $\chi_f(T(G))$ is proven to
be at most $\Delta(G)+2$ for any value of $\Delta(G)$ by Kilakos and Reed  
\cite{KR}. 

Here we prove a topological version, stating that even the strongest of the
above topological bounds is at most $\Delta(G)+2$ for $T(G)$.

\begin{thm}
For any simple graph $G$ the inequality $$\ind(B(T(G)))\le\Delta(G)$$ holds. 
\end{thm}

\gproof
Let $\Delta=\Delta(G)$. 
We prove that $T(G)$ can contain the complete bipartite graph
$K_{2,\Delta+1}$ as a subgraph only if $\Delta\le 3$. 
In the latter case the statement of the theorem follows from the above
mentioned result of Rosenfeld
\cite{Ros} and Vijayaditya \cite{Vij} that verifies 
the original conjecture in this case. For $\Delta>3$ the lack of the above
complete bipartite subgraph proves the theorem by the $K_{\ell,m}$-theorem
of Csorba, Lange, Schurr, and Wa{\ss}mer \cite{CsLSW} quoted in
Section~\ref{top}.  

Assume for a contradiction that $T(G)$ does contain a $K_{2,\Delta+1}$
subgraph, while $\Delta>3$. Let the two sides (color classes) of this
complete bipartite subgraph be denoted by $A$ and $B$, where $|A|=2$ and
$|B|=\Delta+1.$ 

Recall that $V(T(G))=V(G)\cup E(G)$. We will simply denote $V(G)$ by $V$ and
$E(G)$ by $E$ and distinguish among a few cases
according to the size of the intersections of $A$ and $B$ with $V$ and
$E$, respectively.

First observe that $A\cup B\subseteq V(G)$ is impossible, because then the
vertices of $G$ in $A$ should have degree at least $\Delta+1$ contradicting
the definition of $\Delta$. 

So there is some $e\in E$ that belongs to $A\cup B$. 
First assume $e\in B$ and $A\subseteq V$. Then the two vertices in $A$ must be
the two endpoints of $e$. Since there is no other edge both of these vertices
belong to, we must have $|B\cap V|=|B\setminus\{e\}|=\Delta$. Then the
elements of $A$ are adjacent in $G$ (as vertices of $G$) to the $\Delta$
vertices in $|B\cap V|$ plus each other (by $e$), so their degree in $G$ is at
least $\Delta+1$, a contradiction. 

This proves that $A\cap E$ cannot be empty, so we may assume $e\in A$. 
If $|A\cap E|=2$ (that is both elements of $A$ are edges of $G$), say
$A=\{e,f\}\subseteq E$, then $B\cap V\subseteq e\cap f$. If the edges $e$ and
$f$ have no common endpoint, then $B\cap V=\emptyset$ and $\Delta+1=|B|=|B\cap
E|\le 4$, as there are at most $4$ edges that have a common endpoint with both
$e$ and $f$. This contradicts to $\Delta>3$. If, on the other hand, $e$ and
$f$ have a common endpoint $u$, then we still have $|B\cap E|\ge |B|-1=\Delta$
and all these edges except at most one must have $u$ as one of its
endpoints. (One exceptional edge can connect the endpoints of $e$ and $f$
different from $u$.) But then the degree of $u$ is at least $\Delta+1$ in $G$,
a contradiction. 

So we may assume $|A\cap E|=1$. Let $A=\{v,e\}$ where $v\in V$ and $e\in E$. 
If $v$ is an endpoint of $e$ then $|B\cap V|\le 1$ as there are only $2$
vertices $e$ is connected to in $T(G)$ and one of them is $v\notin B$. Thus
$|B\cap E|=|B|-1=\Delta$ and all these $\Delta$ edges have $v$ as one of their
endpoints. But $v$ is also an endpoint of $e\notin B$, so $v$ has degree at
least $\Delta+1$ in $G$, a contradiction again. 
Finally we have to look at the case when $v$ is not an endpoint of $e$. Then
$|B\cap V|\le 2$ since $e$ has only two endpoints and $|B\cap E|\le 2$,
because there are at most $2$ edges containing $v$ as an endpoint and having
the other endpoint at one end of $e$. Thus $\Delta+1=|B|\le 4$ contradicting
the assumption that $\Delta>3$. 
\hfill$\Box$

\begin{rem}
In the proof above we had two cases where the contradiction was
with $\Delta(G)>3$, i.e., where we relied on Rosenfeld's and Vijayaditya's
theorem. The first 
such case is inessential, there we could continue by simply saying that if
$A\cup B\subseteq E$ and the two elements in $A$ are independent edges of $G$,
then getting $|B|=4$ means that the $6$ edges in $A\cup B$ form a $K_4$
subgraph of $G$ which must be a connected component itself and from this point
the 
argument is easy to complete. The second case when we relied on $\Delta>3$ 
is more essential. This is at the very end of the proof and the
$K_{2,4}$ 
produced there can in fact come up in $T(G)$ without forcing the vertices and
edges belonging to it to form a separate component of $G$.
\hfill$\Diamond$
\end{rem}

\section{On Hedetniemi's conjecture} \label{secHedet}

For two graphs $F$ and $G$ their direct (or categorical) product
$F\times G$ is defined on vertex set $V(F)\times V(G)$ such that two vertices
$(f_1,g_1)$ and $(f_2,g_2)$ are adjacent if and only if $\{f_1,f_2\}\in E(F)$
and $\{g_1,g_2\}\in E(G)$. Let $F$ and $G$ be simple graphs. It is easy to
check that $\chi(F\times G)\le \min\{\chi(F),\chi(G)\}$ (simply color vertex
$(f,g)$ of $F\times G$ with the 
color of $f$ to obtain a proper coloring with $\chi(F)$ colors), and
Hedetniemi's conjecture states that equality holds. This conjecture is wide
open, the 
major special case proven is when the right hand side is $4$ \cite{EZS}. In
fact even that is not known whether the function $f(t):=\min\{\chi(F\times G):
\chi(F)\ge t, \chi(G)\ge t\}$ goes to infinity with $t$ or not. (Though rather
surprisingly, if it does not, then it remains below $10$.) For further
information and references we refer the reader to the excellent 
recent survey by Tardif \cite{Tardsurv}. 

Concerning relaxations involving the fractional chromatic number, Tardif
proved in \cite{THed1} that $\chi(F\times G)\ge {1\over
  2}\min\{\chi_f(F),\chi_f(G)\}$ and in \cite{THed2} that 
$\chi_f(F\times G)\ge {1\over
  4}\min\{\chi_f(F),\chi_f(G)\}$, where $\chi_f$ stands for the fractional
chromatic number. Note that while the first of these inequalities is a
weakening of Hedetniemi's conjecture, the second is only an analog, although
if the exact value of the multiplicative constant is ignored then it also
implies the first one.  

Another relaxation mentioned in Tardif's survey \cite{Tardsurv} is due to Hell
\cite{Hell}. It already connects Hedetniemi's conjecture to Lov\'asz's
topological 
lower bound on the chromatic number. In particular, in Tardif's formulation,
Hell shows that if $F$ and $G$ are two graphs for which Lov\'asz's bound is
tight then $\chi(F\times G)= \min\{\chi(F),\chi(G)\}$.

Along these lines we state the following topological analog of Hedetniemi's
conjecture.  

\begin{thm} \label{hedrel}
$$\coind(B(F\times G))=\min\{\coind(B(F)),\coind(B(G))\}.$$
\end{thm} 
 
\gproof
It is true for any pair of graphs $F$ and $G$ that a homomorphism from
$F\times G$ exists both to $F$ and $G$ (simply by taking
projections). Assume $\coind(B(F\times G))=h$. Then, by Theorem B, there is some
Borsuk graph 
$B(h+1,\alpha)$ which homomorphically maps into $F\times G$. Combining this
homomorphism with either of the projection homomorphisms mentioned above, we
get a homomorphism from $B(h+1,\alpha)$ to $F$ and to $G$, respectively. Thus
$\coind(B(F))\ge h$ and $\coind(B(G))\ge h$ also holds. This proves
$\coind(B(F\times G))\le\min\{\coind(B(F)),\coind(B(G))\}.$ 

To prove the reverse inequality let
$d=\min\{\coind(B(F)),\coind(B(G))\}$. Then by Theorem B there is some large
enough $\alpha<2$ for which $B(d+1,\alpha)$ admits a homomorphism $f$ to $F$
and also a homomorphism $g$ to $G$. But then the function which maps every
vertex $v\in V(B(d+1,\alpha))$ to $(f(v),g(v))$ is a homomorphism of
$B(d+1,\alpha)$ to $F\times G$ and its existence proves $\coind(F\times G)\ge
d$. Therefore we have $\coind(B(F\times
G))\ge\min\{\coind(B(F)),\coind(B(G))\}$ completing the proof.  
\hfill$\Box$

\begin{rem}
P\'eter Csorba \cite{CsPprivat} observed that a result of Kozlov 
(equation 2.4.2 in \cite{Koz}; cf. also Hell \cite{Hell}) combined with results
from Csorba \cite{Csorbhom} implies 
that $B(F\times G)$ is $\2$-homotopy equivalent to $B(F)\times B(G)$ for every
pair of graphs $F$ and $G$. The product $A\times
B$ here is meant to be the product topological space equipped with the $\2$-map
$\nu: (x,y)\mapsto (\nu_A(x),\nu_B(y))$ where $\nu_A$ and $\nu_B$ are  the
respective $\2$-maps on the $\2$-spaces $A$ and $B$.   
From this observation an alternative proof for Theorem~\ref{hedrel} can easily
be obtained.   
\hfill$\Diamond$
\end{rem}

\section{On Hadwiger's conjecture} \label{secHad}

Hadwiger's conjecture states that if $\chi(G)\ge t$, then $G$ contains a
$K_t$ minor. It is essentially trivial for $t\le 3$,
relatively easy to prove for $t=4$, 
known to be equivalent to the Four Color Theorem for $t=5$, and proven
to be so also for $t=6$ \cite{RST}; see \cite{Toft} for an excellent 
survey. For $t\ge 7$ the conjecture is open and is widely considered as one of
the most important open problems in graph theory. Even a linear approximation
is not proven, that is, it is not known whether there exists a constant $c$
such that $\chi(G)\ge t$ implies that $G$ contains a complete
minor on $ct$ vertices. Such a result is proven for the fractional chromatic
number in place of the chromatic number with $c=1/2$ by Reed and Seymour
\cite{RS}. Stating in the counterpositive form they proved that if a graph
contains no complete minor on $m+1$ vertices then its fractional chromatic
number is at most $2m$. In \cite{ST} it was observed that an analogous
statement immediately follows from the $K_{\ell,m}$-theorem for the
topological lower bounds on the chromatic number. Namely, if $G$
contains no $K_{m+1}$ minor, then $\ind(B(G))+2< 2m$. 

Although Hadwiger's conjecture is wide open, a strengthening, called the ``odd
Hadwiger conjecture'' received much attention in recent years, see
\cite{GGRSV, KenK, KKR, KKS}. To state it we need the concept of an
odd $K_m$ minor.   

\begin{defi}
An odd $K_m$ minor of graph $G$ is formed by $m$ vertex disjoint trees
$T_1, \dots, T_m$ in $G$ that have the following additional properties. 

The vertices in these trees can be simultaneously $2$-colored such that 

1. Each tree $T_i$ is properly colored;

2. For every pair $i\ne j$, $1\le i,j\le m$,
there is an edge between a vertex of $T_i$ and a vertex of $T_j$ that have the
same color.
\end{defi}

The odd Hadwiger conjecture was suggested by Gerards and Seymour
(cf.\ \cite{JT}) and it states that if $\chi(G)\ge t$ then $G$ must contain an
odd $K_t$ minor. 

In some cases the known results about this conjecture show surprising
similarities with those known about Hadwiger's original conjecture. In
particular, Kawarabayashi and Reed \cite{KKR} have proved an analog of the
Reed-Seymour theorem, namely, they showed that if $G$ does not
contain an odd $K_m$ minor then the fractional chromatic number of $G$ is at
most $2m$. This suggests the question whether one can prove that graphs with
no odd $K_m$ minor satisfy $\ind(B(G))+2\le 2m$. It is clear that now we cannot
get this just from the $K_{\ell,m}$-theorem, since its conclusion holds for
large complete bipartite graphs that contain no large odd complete minors. 
Though we also did not succeed to get something similar from the Zig-zag
Theorem, we wonder whether its conclusion would already imply such a
statement.  

\begin{tquestion}
Is there some constant $c$ for which the following is
true? If every proper coloring of a graph $G$ satisfies the
conclusion of the Zig-zag Theorem (with parameter $t$) then $G$ contains an
odd complete minor on $ct$ vertices.
\end{tquestion}

In the following section we prove that the odd Hadwiger conjecture holds for
some graph families that have their chromatic number equal to its topological
lower bound, while the fractional chromatic number is much smaller.

\section{The odd Hadwiger conjecture for large enough Kneser graphs and
  generalized Mycielski graphs} \label{secoddHad}

Recall that the Kneser graph $\KG(n,k)$ is defined for $n>2k$ on all
$k$-subsets of an $n$-element set as vertices, where two of these form an edge
iff they are disjoint. The chromatic number of graph $\KG(n,k)$ is
$n-2k+2$ \cite{LLKn}, while their fractional chromatic number is only $n/k$,
cf. \cite{EKR, SchU}. 
We are going to prove the following result.

\begin{thm} \label{knh0}
If $t=n-2k+2$ is fixed and $n$ is large enough, then the $t$-chromatic Kneser
graph $\KG(n,k)$ contains an odd $K_t$ minor. 
\end{thm}

Recall that a topological $K_r$ subgraph in a graph $G$ is a collection of $r$
{\it branching vertices} together with $r\choose 2$ vertex-disjoint paths in
$G$ connecting all pairs of the branching vertices. We call a topological
$K_r$ subgraph {\it odd} if all the latter $r\choose 2$ paths are odd, i.e.,
they contain an odd number of edges. A famous conjecture stronger than
Hadwiger's was due to Haj\'os claiming that every graph of chromatic number
$t$ contains a topological $K_t$ subgraph. This was disproved by Catlin
\cite{Cat} for $t\ge 7$, cf.\ also \cite{Toft}. (It is known to hold for $t\le
4$ when it is actually 
equivalent to Hadwiger's conjecture, and is still open for $t=5, 6$.) Several
other counterexamples can be found in a more recent paper by Thomassen
\cite{ThomHaj}.

Since every odd topological $K_r$ subgraph gives rise to an odd $K_r$ minor,
and   
since the odd Hadwiger conjecture is trivial if the chromatic number is less
than $4$ and is also known to hold when it is equal to $4$ (the latter was
proven by Catlin \cite{Cat}), Theorem~\ref{knh0} immediately follows from the
following result. 

\begin{thm} \label{knoh}
If $t=n-2k+2\ge 5$ is fixed and $n$ is large enough then the Kneser graph
$\KG(n,k)$ contains an odd topological $K_t$ subgraph. 
\end{thm}

\gproof
Arrange the $n$ points $1, 2,\dots, n$ on a circle and let their $k$-subsets
be identified with the vertices of $\KG(n,k)$. 

A $k$-subset formed by $k$ cyclically consecutive points on the circle will be
called a {\it short arc}, while a {\it long arc} is formed by a set of
$\ell:=\left\lfloor 
  {{n-1}\over 2}\right\rfloor$ cyclically consecutive points. The first 
point of a long arc is meant to be its first element when the arc is traversed
in a clockwise order along the circle. The relative position of a $k$-subset of
a long arc within the long arc will be called its {\it pattern} if it contains
the first element of the long arc. 
Thus for example, if $k=3$, then the subset $\{1,3,7\}$ has the
same pattern in the long arc starting with $1$ as the subset $\{n-1,1,5\}$
in the long arc starting with $n-1$. Note that a pattern in a given long
arc defines a vertex of $\KG(n,k)$, and if we have two different pairs, both 
consisting of a long arc and a pattern, then these define distinct vertices of
$\KG(n,k)$. This is ensured by the condition that the first element of the
long arc is always in the $k$-subset defined by a pattern. Note also that for
such vertices it is meaningful to speak about the {\it pattern of the
  vertex}, since a $k$-subset that fits into a long arc defines the long arc
with which it has the same starting vertex, i.e., the one in which it has a
pattern.   

Select $t=n-2k+2$ short arcs (there are $n$ altogether) and fix them. These
will 
be the branching vertices of our odd topological $K_t$. Call a pattern
\emph{good} if it is not identical with the first $k$ vertices in the long arc. 
(In other words, vertices with a good pattern are not short arcs.) 
Next we select a different good pattern for each pair of the branching vertices. 
First we will show that this is possible, 
and then we show that between any two branching
vertices there is a path of odd length, all inner vertices of which have the
same pattern, namely the one attached to the given pair of vertices. These
paths 
will then be automatically disjoint as their inner
vertices have different patterns. They also cannot touch other branching
vertices than their endpoints, since the branching vertices form short arcs
and they are excluded from the set of good patterns.  
 
The number of good patterns is easily seen to be
${\left\lfloor{(n-3)/2}\right\rfloor\choose {k-1}}-1$. 
This is equal to ${{(n-3)/2}\choose{(t-3)/2}}-1$ if $n$ is odd and  
to ${{(n-4)/2}\choose{(t-4)/2}}-1$ if $n$ is even. 

We can select a different good pattern for all pairs of branching vertices if 
the above expression is not less than $t\choose 2$. Since $t$ is fixed, the
latter number is constant, while the above expressions go to infinity with $n$
whenever $t\ge 5$. This proves that for large enough $n$ the required
inequality holds. 

It remains to prove that between any two branching vertices there exists a
path of odd length with all inner nodes having an arbitrarily fixed good
pattern. To this end, first observe that for any two given long arcs, $a$ and
$b$, one can find a sequence of long arcs $a=s_0, s_1, \dots, s_r=b$, such
that $s_i\cap s_{i+1}=\emptyset$ for all $i=0,1,\dots,r-1$ and $r$ is
even. In other words, there is an even length path between any two vertices in
$\KG(n,\ell)[L]$, where $\KG(n,\ell)[L]$ is the subgraph of $\KG(n,\ell)$
induced by the set $L$ of vertices that form long arcs. 
This statement is true because two closest long arcs, i.e., two long
arcs with symmetric difference $2$ still have a long arc in the complement of
their union. Thus with two steps (a step meaning going from one vertex of
$\KG(n,\ell)[L]$ to another along an edge) we can shift any long arc along our
circle by $1$. Therefore we can realize any shift with an even number of
steps.  
Given two branching points, i.e., two short arcs $x$ and $y$, choose $a$ to be
a long arc disjoint from $x$ and $b$ to be a long arc containing $y$. Consider
the above sequence of long arcs between $a$ and $b$ and then substitute each
long arc of the sequence by the vertex of $\KG(n,k)$ contained in the given
long arc and having the pattern attached to the pair of branching points
$(x,y)$ (while $b$ is substituted by $y$). Adding $x$ to the beginning of this
sequence we obtain the required odd
length path in $\KG(n,k)$ between vertices $x$ and $y$ completing the proof. 
\hfill$\Box$

Schrijver \cite{Schr} defined a beautiful family of graphs, that appear as
induced subgraphs of Kneser graphs and share some of their important
properties. 

\begin{defi} {\rm (\cite{Schr})}
The Schrijver graph $\SG(n,k)$ is defined as the induced subgraph of the
Kneser graph $\KG(n,k)$ on the vertices 
$$
  V(SG(n,k)) = 
    \big\{a\subseteq {[n]\choose k}: \forall i\ \{i,i+1\}\nsubseteq a,
    \{1,n\}\nsubseteq a\big\}.$$
\end{defi}
 
It is proven in \cite{Schr} that the chromatic number of $\SG(n,k)$ is also
$n-2k+2$ as for $\KG(n,k)$, moreover, $\SG(n,k)$ is vertex color-critical. 
Talbot \cite{Talb} determined the independence number of $\SG(n,k)$ which
easily implies that $\chi_f(\SG(n,k))=\chi_f(KG(n,k))$, too. 

It is easy to see, that if $n$ is odd, then choosing the cyclic permutation on
our circle at
the beginning of the proof of Theorem~\ref{knoh} as 
$ 1, 3, 5, \dots, n, 2, 4, \dots, n - 1 $, each
long arc will be such that neither a set $\{i,i+1\}$, nor $\{n,1\}$ will be
contained in it. Thus any $k$-subset of any long arc will be a vertex of
$\SG(n,k)$ and the proof goes through for $\SG(n,k)$ just as it did for
$\KG(n,k)$. In case $n$ is even, we
can simply ignore the point $n$ and fix the circle as above on the elements
$1,\dots,n-1$ only. Observing that the proof would allow more than $t$
branching points, too, we apply the above argument for $t+1$ branching points
that goes again through the same way. Thus we obtain the following
strengthening of Theorems~\ref{knh0} and \ref{knoh}.

\begin{cor}
If $t=n-2k+2\ge 5$ is fixed and $n$ is large enough, then the $t$-chromatic
Schrijver graph $\SG(n,k)$ contains an odd topological $K_t$ subgraph. In
particular, for any fixed $t=n-2k+2$ and $n$ large enough 
$\SG(n,k)$ contains an odd $K_t$ minor.
\end{cor}

Generalized Mycielski graphs form another family of graphs where
topological lower bounds on the chromatic number give sharp estimates, while
the fractional chromatic number is far below the chromatic number
\cite{Tar}. 

\medskip
The $r$-level generalized Mycielskian $M_r(G)$ of a graph $G$ is defined on
the vertex set $$V(G)\times \{0,1,\dots, r-1\}\cup \{z\}$$ with edge set
$$\displaylines{\quad
E(M_r(G))=\bigl\{\{(u,i),(v,j)\}: \{u,v\}\in E(G)\textup{ and } 
(|i-j|=1\textup{ or }i=j=0)\bigr\}\cup{}
\hfill\cr\hfill
\bigl\{\{(u,r-1),z\}: u\in V(G)\bigr\}.
\quad}$$
Mycielski graphs are
usually meant to be the graphs obtained from $K_2$ by an iterative use of the
above general Mycielski construction with $r=2$.  

Results of Stiebitz \cite{Stieb} (cf.\ also \cite{GyJS, Mat}) generalized by
Csorba \cite{CsorbMyc} imply that if the box complex $B(G)$ of a graph $G$ is
homotopy equivalent to a sphere $\S^h$ (this is the case for complete graphs
and more generally for all Schrijver graphs, see \cite{BjdeL}), then for any
positive 
integer $r$, the box complex $B(M_r(G))$ is homotopy equivalent to $\S^{h+1}$,
therefore $\ind(B(M_r(G)))=\ind(B(G))+1$ holds. 
In particular, if the above homotopy equivalence holds 
and the topological lower bound (in this case the four lower bounds we
discussed coincide) of the chromatic number is tight (this also happens for
all Schrijver graphs), then it is $1$ more and
also tight for $M_r(G)$. (Note that there are graphs with
$\chi(M_r(G))=\chi(G)$, an example given in \cite{Tar} is the complement of
the $7$-cycle with $r=3$. Another example is given in \cite{CsorbMyc}.) 

Concerning the odd Hadwiger conjecture we prove the following. 

\begin{prop}
If $G$ contains an odd $K_t$ minor then $M_r(G)$ contains an odd $K_{t+1}$
minor for every $r\ge 1.$ 
\end{prop}

\gproof
We may assume that $G$ is connected and that $ r\ge 2 $. 
Consider $G$ as the subgraph induced on vertices $(v,0)$ of $M_r(G)$ and the
$t$ vertex disjoint trees $T_1, \dots, T_t$ with their $2$-coloring that give
an odd $K_t$ minor in this induced subgraph $G$. Notice that if some of these
$t$ trees have only one vertex then they are all colored the same, say blue. 

Now take an arbitrary spanning tree $T_{t+1}$ on the vertices in the set
$\{(v,i): i>0\}\cup \{z\}$ and its proper $2$-coloring that gives color blue
to all vertices of the 
form $(v,1)$. (Such a coloring is valid as the vertices $\{(v,i): i>0\}\cup
\{z\}$ induce a bipartite subgraph in $M_r(G)$ in which the distance between
any two vertices $\{v,1\}, \{v',1\}$ is even.) 
It remains to show only that all trees $T_i$ with $i\le t$ have a
blue colored vertex that has a neighbor among the vertices $(v,1)$. But this
is almost obvious: By the connectedness of $G$ every vertex $(u,0)$ has some
neighbor of the form $(v,1)$ and all the trees $T_1,\dots,T_t$ either have an
edge and then one of its endpoints is necessarily blue or it is a one-point
tree, but then it is blue by the above observation. So we have an
odd $K_{t+1}$ minor. \hfill$\Box$

\bigskip
\bigskip
\par\noindent
{\bf Acknowledgement.} 

\par\noindent
An inspiring conversation with Claude Tardif and G\'abor Tardos is
gratefully acknowledged. We also thank for a very useful conversation to
P\'eter Csorba.

\end{document}